\documentclass[12pt]{article}

\usepackage{a4}
\usepackage{amsthm}
\usepackage{amsfonts}
\usepackage{epsfig,ifpdf}
\newtheorem{theorem}{Theorem}
\newtheorem{lemma}[theorem]{Lemma}
\ifpdf
\fi

\newcommand{\G}{{\cal G}}

\begin{document}
\title{Long cycles in fullerene graphs}
\author{%
Daniel Kr{\'a}l'$^a$\footnote{Institute for Theoretical Computer Science ({\sc iti})
is supported by Ministry of Education of the Czech Republic as project 1M0545.}
\and
Ond{\v r}ej Pangr{\'a}c$^a$
\and
Jean-S{\'e}bastien Sereni$^a$\footnote{This author is supported by the
European project {\sc ist fet Aeolus.}}
\and
Riste {\v S}krekovski$^b$\footnote{Supported in part by
Ministry of Science and Technology of Slovenia, Research Program
P1-0297.}}
\date{}
\maketitle
\begin{center}
$^a$ Department of Applied Mathematics ({\sc kam}) and
Institute for Theoretical Computer Science ({\sc iti}),
Faculty of Mathematics and Physics, Charles University, Malostransk\'e N\'am\v est\'i 25,
118 00 Prague, Czech Republic.\\
E-mails: \{{\tt kral,pangrac,sereni}\}{\tt @kam.mff.cuni.cz}.\\
$^b$ Department of Mathematics, University of Ljubljana,
Jedranska 19, 1111 Ljubljana, Slovenia.\\
\end{center}

\begin{abstract}
It is conjectured that every fullerene graph is hamiltonian.
Jendrol' and Owens proved [J. Math. Chem. 18 (1995), pp.~83--90]
that every fullerene graph on $n$ vertices has a cycle
of length at least $4n/5$.
In this paper, we improve this bound to $5n/6-2/3$.
\end{abstract}

%
%

\section{Introduction}

\emph{Fullerenes} are carbon-cage molecules comprised of carbon atoms
that are arranged on a sphere with twelve pentagon-faces and other
hexagon-faces. The~ico\-sahedral $C_{60}$, well known as Buckministerfullerene,
was found by Kroto et~al.~\cite{KHBCS}, and later confirmed through experiments by Kr\"{a}tchmer et~al.~\cite{KLFH} and Taylor et al.~\cite{THAK}.
Since the discovery of the first fullerene molecule, the
fullerenes have been objects of interest to scientists in many disciplines.

Many properties of fullerene molecules can be studied using mathematical tools
and results. Thus, \emph{fullerene graphs} were defined as cubic (i.e.~$3$-regular) planar 3-connected graphs with pentagonal and hexagonal faces. Such graphs are suitable models for fullerene molecules: carbon atoms are represented by vertices of the graph, whereas the edges represent bonds between adjacent atoms.
It is known that there exists a fullerene graph on $n$ vertices for every
even $n\ge 20$, $n\not=22$. See the monograph of Fowler and Manolpoulos~\cite{FM} for more information on fullerenes.

The hamiltonicity of planar 3-connected cubic graph has been attracting much
interest of mathematicians since Tait~\cite{Tait} in 1878 gave a short and elegant
(but also false) proof of the Four Color Theorem based on the "fact" that planar 3-connected
cubic graphs are hamiltonian. The missing detail of the proof was precisely
the previously mentioned ``fact'' which became known as Tait's Conjecture. Later, Tutte~\cite{Tutte} disproved Tait's Conjecture.

The hamiltonicity of various subclasses of 3-connected planar cubic graphs was
additionally investigated. Gr\"unbaum and Zaks~\cite{GZ} asked whether
the graphs in the family $\G_3(p, q)$  of 3-connected cubic planar graphs
whose faces are of size $p$ and $q$ with $p<q$ are hamiltonian for any $p,q$.
Note that $p \in \{3,4,5\}$ by Euler's formula.
Also note that fullerene graphs correspond to $\G_3(5,6)$.  Goodey~\cite{Goodey1,Goodey2} has proved that all graphs contained in
$\G_3(3,6)$ and $\G_3(4,6)$ are hamiltonian. Zaks~\cite{Zaks} found non-hamiltonian graphs in the family $\G_3(5, k)$ for $k\ge  7$.
Similarly, Walther~\cite{Walther} showed that families $\G_3(3, q)$
for $7\le q\le 10$ and $\G_3(4, 2k + 1)$ for $k\ge 3$, contain
non-hamiltonian graphs.
For more results in this area, also see~\cite{O1,O2,O3, Tkac}.

Let us restrict our attention to $\G_3(5,6)$. Ewald~\cite{Ewald}
proved that every fullerene graph contains a cycle which meets
every face of $G$. This implies that there is a cycle through at
least $n/3$ of the vertices of any fullerene graph on $n$
vertices. It is well known that each graph $G\in\G_3(5,6)$ has a
dominating cycle $C$, i.e. a cycle $C$ such that each edge of $G$
has an end-vertex on $C$ (in particular, a Tutte cycle is dominating since $G$ is
cyclically $4$-edge-connected).
This immediately improves the
bound from $n/3$ to $3n/4$. Jendrol' and Owens~\cite{JO} gave a
better bound of $4n/5$. In this paper, we improve the bound to
$5n/6-2/3$.

%
%

\section{Preliminary observations}

We follow the terminology of Jendrol' and Owens~\cite{JO}. We consider a longest
cycle $C$ of a fullerene graph; a vertex contained in $C$ is
{\em black} and a vertex not contained in $C$ is {\em white}.
Our aim is to show that there are at most $n/6+2/3$ white vertices
for an $n$-vertex fullerene graph $G$. The following was
shown~\cite{JO}.

\begin{lemma}
\label{lm-path}
Let $G$ be a fullerene graph and $C$ a longest cycle in $G$.
The graph $G$ contains no path comprised of three white vertices.
\end{lemma}

\begin{figure}
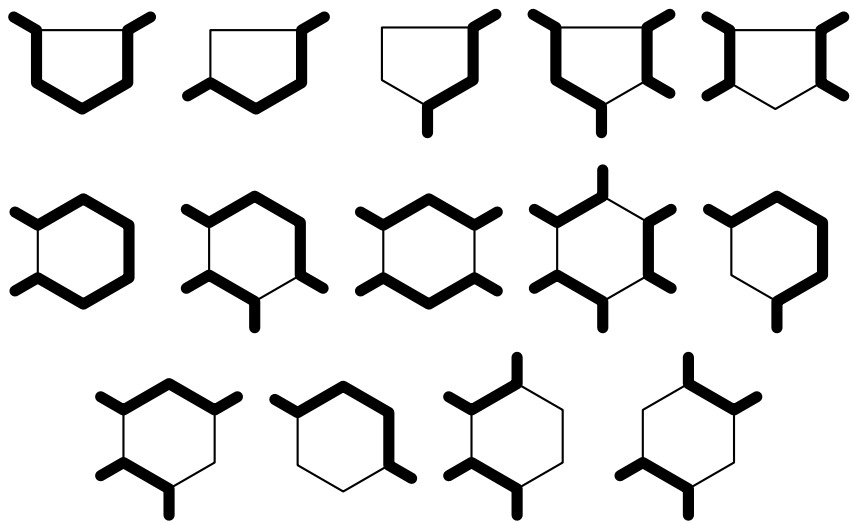

\begin{center}
\epsfbox{ful-ham.10}\hspace{2mm}
\epsfbox{ful-ham.11}\hspace{2mm}
\epsfbox{ful-ham.12}\hspace{2mm}
\epsfbox{ful-ham.13}\hspace{2mm}
\epsfbox{ful-ham.14}\\[3mm]
\epsfbox{ful-ham.1}\hspace{2mm}
\epsfbox{ful-ham.2}\hspace{2mm}
\epsfbox{ful-ham.3}\hspace{2mm}
\epsfbox{ful-ham.4}\hspace{2mm}
\epsfbox{ful-ham.5}\\[2mm]
\epsfbox{ful-ham.6}\hspace{2mm}
\epsfbox{ful-ham.7}\hspace{2mm}
\epsfbox{ful-ham.8}\hspace{2mm}
\epsfbox{ful-ham.9}
\end{center}
\caption{The possible ways for a cycle $C$ to traverse a face of size five
or six (up to symmetry) without forming a path of three white vertices.
The cycle $C$ is indicated by bold edges.}
\label{fig-possible}
\end{figure}

Lemma~\ref{lm-path} implies that no face of $G$ is incident
with more than two white vertices (see Figure~\ref{fig-possible}
for all possibilities, up to symmetry, how the cycle $C$ can
traverse a face of $G$).
The faces incident with two white vertices are
called {\em white} and the faces incident with no white vertices
are called {\em black}. Let us now observe the following simple
fact.

\begin{lemma}
\label{lm-white5}
If $C$ is a longest cycle in a fullerene graph $G$, then there
are no white faces of size five.
\end{lemma}

\begin{figure}
\begin{center}
\epsfbox{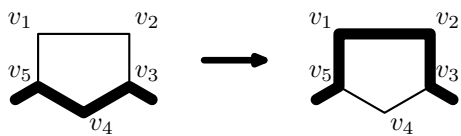}
\end{center}
\caption{Prolonging the cycle $C$ if the graph $G$ contains a white face
of size five.}
\label{fig-white5}
\end{figure}

\begin{proof}
Assume that there is a white face $v_1v_2v_3v_4v_5$. By symmetry
we can assume, that the cycle $C$ contains the path $v_3v_4v_5$.
Replacing the path $v_3v_4v_5$ with the path $v_3v_2v_1v_5$
(see Figure~\ref{fig-white5}) yields a cycle of $G$ longer than $C$,
a contradiction.
\end{proof}

In the sequel, we use the following notion.
Given a face $f$ with vertices $v_1,v_2,\ldots,v_k$ (in cyclic order), we
let $f_{i,i+1}$ be the face different from $f$ that contains the edge
$v_{i}v_{i+1}$ (the indices are taken modulo $k$).

%
%

\section{Initial charge and discharging rules}
\label{sect-rules}

Using a discharging argument, we argue that the number of white vertices
with respect to a longest cycle $C$ in an $n$-vertex fullerene graph $G$
is at most $n/6+2/3$. Fix such a cycle $C$. Each white vertex
initially receives $3$ units of charge. Next, each white vertex
sends $1$ unit of charge to each of its three incident faces.
Observe that each white face has $2$ units of charge, each black face
has no charge and each remaining face has $1$ unit of charge each.

The charge is now redistributed based on the following rules (the indices
are taken modulo the length of the considered face where appropriate).

\begin{description}
\item[Rule A] A black face $f_0=v_1\ldots v_6$ receives $1/2$ unit
of charge from the face $f_{i,i+1}$ if the path $v_{i-1}v_iv_{i+1}v_{i+2}$
is contained in the cycle $C$ and the face $f$ is white.
\item[Rule B] A black face $f_0=v_1\ldots v_6$ receives $1$ unit
of charge from the face $f_{i,i+1}$ if the edge $v_iv_{i+1}$ is contained
in the cycle $C$, neither the edge $v_{i-1}v_i$ nor the edge $v_{i+1}v_{i+2}$
is contained in $C$ and the face $f$ is white.
\end{description}

\begin{figure}
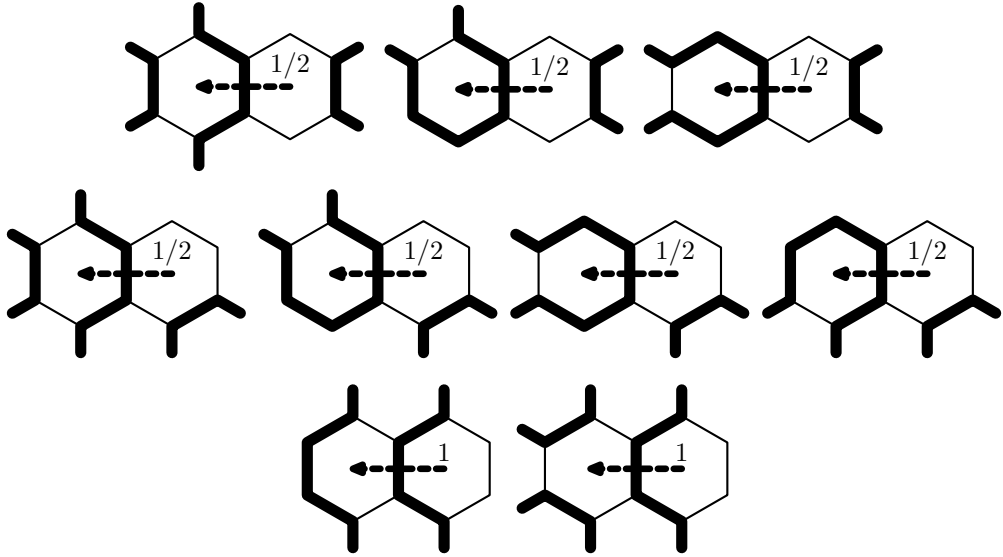

\begin{center}
\epsfbox{ful-ham.16}\hspace{1mm}
\epsfbox{ful-ham.17}\hspace{1mm}
\epsfbox{ful-ham.18}\\[2mm]
\epsfbox{ful-ham.19}
\epsfbox{ful-ham.20}
\epsfbox{ful-ham.21}
\epsfbox{ful-ham.22}\\[3mm]
\epsfbox{ful-ham.23}\hspace{2mm}
\epsfbox{ful-ham.24}
\end{center}
\caption{Configurations (up to symmetry) to which Rules A and B are applied.}
\label{fig-rules}
\end{figure}

\noindent The Rules A and B are illustrated in Figure~\ref{fig-rules}.

In Sections~\ref{sect-white} and~\ref{sect-black}, we show that each
face has at most $1$ unit of charge after applying Rules A and B.
Based on this fact, we conclude in Section~\ref{sect-main} that the number
of white vertices is at most $f/3$ where $f$ is the number of faces of $G$.
The bound on the length of the cycle $C$ will then follow.

%
%

\section{Final charge of white faces}
\label{sect-white}

In this section, we analyze the final amount of charge of white faces.
By Lemma~\ref{lm-white5}, we can restrict our attention to faces of
size six.

\begin{lemma}
\label{lm-white6-para}
Let $C$ be a longest cycle of a fullerene graph $G$. Assume that
the discharging rules as described in Section~\ref{sect-rules}
have been applied. If $f=v_1v_2v_3v_4v_5v_6$ is a white face of $G$
such that the edges $v_2v_3$ and $v_5v_6$ are contained in $C$,
then the final amount of charge of $f$ is $1$ unit.
\end{lemma}

\begin{figure}
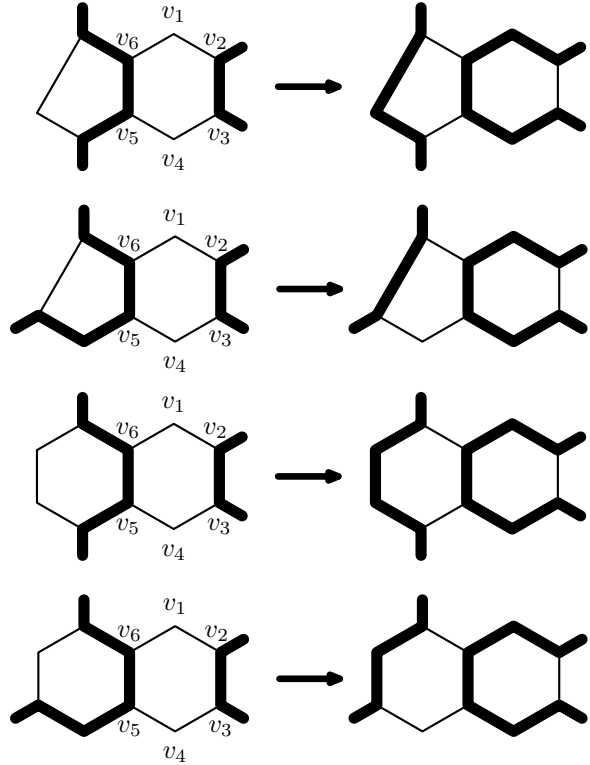

\begin{center}
\epsfbox{ful-ham.27}\\[4mm]
\epsfbox{ful-ham.28}\\[2mm]
\epsfbox{ful-ham.25}\\[4mm]
\epsfbox{ful-ham.26}
\end{center}
\caption{Configurations analyzed in the proof of Lemma~\ref{lm-white6-para}.}
\label{fig-white6-para}
\end{figure}

\begin{proof}
The initial amount of charge of the face $f$ is $2$ units.
If both the faces $f_{23}$ and $f_{56}$ are black faces of size six,
then the face $f$
sends $1/2$ unit of charge to each of them by Rule A and thus
its final amount of charge is $1$ unit.

Assume that the face $f_{56}$ is not a black face of size six.
Hence, the graph $G$,
up to symmetry, contains one of the configurations depicted
in the left column of Figure~\ref{fig-white6-para}.
Rerouting the cycle $C$ as indicated
in the figure yields a cycle longer than $C$, a contradiction.
Since our arguments translate to the case where
the face $f_{23}$ is not a black face of size six,
the proof of the lemma is finished.
\end{proof}

\begin{lemma}
\label{lm-white6-orto1}
Let $C$ be a longest cycle of a fullerene graph $G$. Assume that
the discharging rules as described in Section~\ref{sect-rules}
have been applied. If $f=v_1v_2v_3v_4v_5v_6$ is a white face of $G$
such that the edges $v_4v_5$, $v_5v_6$ and $v_6v_1$ are contained in $C$,
then the final amount of charge of $f$ is $1$ unit.
\end{lemma}

\begin{figure}
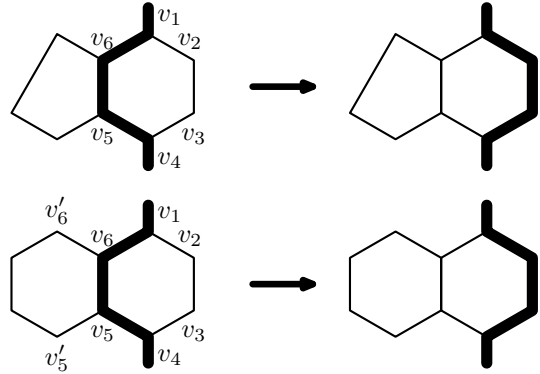

\begin{center}
\epsfbox{ful-ham.36}\\[3mm]
\epsfbox{ful-ham.35}
\end{center}
\caption{Configurations analyzed in the proof of Lemma~\ref{lm-white6-orto1}.}
\label{fig-white6-orto1}
\end{figure}

\begin{proof}
First, suppose that the face $f_{56}$ has size five.
Rerouting the cycle $C$ as indicated in the top line of
Figure~\ref{fig-white6-orto1} yields a face of size five incident
with two or more white vertices (the vertices $v_5$ and $v_6$ become
white). This is excluded by Lemma~\ref{lm-white5}.

We conclude that the face $f_{56}$ has size six.
For $i\in\{5,6\}$, let $v'_i$ be the neighbor of the vertex $v_i$ that
is not incident with the face $f$.
The vertex $v'_5$ cannot be white: otherwise, rerouting
the cycle $C$ as indicated in the bottom line of Figure~\ref{fig-white6-orto1}
yields a path formed by three white vertices. This is impossible
by Lemma~\ref{lm-path}. Thus, the vertex $v'_5$ is black.
Similarly, the vertex $v'_6$ is black. Consequently, the face
$f_{56}$ is black and by Rule B, the face $f_{56}$ receives
$1$ unit of charge from the face $f$. Since the face $f$
sends charge to no other face and its initial amount of charge
is $2$ units, its final amount of charge is $1$ unit.
\end{proof}

\begin{lemma}
\label{lm-white6-orto2}
Let $C$ be a longest cycle of a fullerene graph $G$. Assume that
the discharging rules as described in Section~\ref{sect-rules}
have been applied. If $f=v_1v_2v_3v_4v_5v_6$ is a white face of $G$
such that the edges $v_3v_4$ and $v_5v_6$ are contained in $C$ and
the edge $v_4v_5$ is not contained in $C$,
then the final amount of charge of $f$ is $1$ unit.
\end{lemma}

\begin{figure}
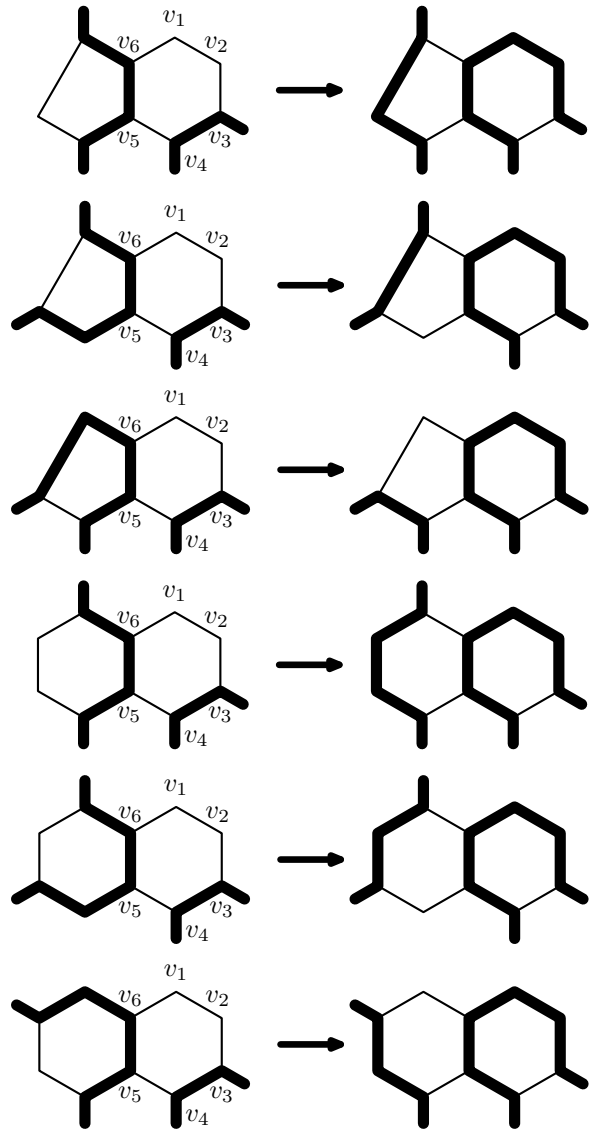

\begin{center}
\epsfbox{ful-ham.32}\\[3mm]
\epsfbox{ful-ham.33}\\[2mm]
\epsfbox{ful-ham.34}\\[3mm]
\epsfbox{ful-ham.29}\\[3mm]
\epsfbox{ful-ham.30}\\[2mm]
\epsfbox{ful-ham.31}
\end{center}
\caption{Configurations analyzed in the proof of Lemma~\ref{lm-white6-orto2}.}
\label{fig-white6-orto2}
\end{figure}

\begin{proof}
The initial amount of charge of the face $f$ is $2$ units.
If both the face $f_{34}$ and $f_{56}$ are black faces of size six,
then the face $f$
sends $1/2$ unit of charge to each of them by Rule A and thus
its final amount of charge is $1$ unit.

Assume that the face $f_{56}$ is not a black face of size six.
Hence, the graph $G$,
up to symmetry, contains one of the configurations depicted
in the left column of Figure~\ref{fig-white6-orto2}.
Rerouting the cycle $C$ as indicated
in the figure yields a cycle of $G$ longer than $C$, a contradiction.
Since our arguments translate to the case
where the face $f_{34}$ is not a black face of size six,
the proof of the lemma is finished.
\end{proof}

Lemmas~\ref{lm-white5}, \ref{lm-white6-para},
\ref{lm-white6-orto1} and~\ref{lm-white6-orto2}
yield the following.

\begin{lemma}
\label{lm-white}
Let $C$ be a longest cycle of a fullerene graph $G$. Assume that
the discharging rules as described in Section~\ref{sect-rules}
have been applied. The final amount of charge of any white face
of $G$ is $1$ unit.
\end{lemma}

%
%

\section{Final charge of black faces}
\label{sect-black}

\begin{figure}
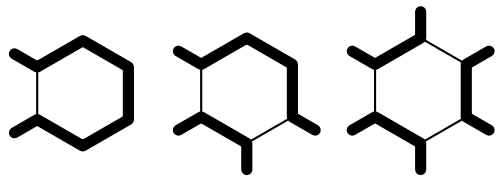

\begin{center}
\epsfbox{ful-ham.1}\hspace{1mm}
\epsfbox{ful-ham.2}\hspace{2mm}
\epsfbox{ful-ham.4}
\end{center}
\caption{Black faces of size six that could receive more than $1$ unit of charge.}
\label{fig-black}
\end{figure}

This section is devoted to the analysis of the final charge of black
faces. Since no black face of size five receives any charge,
we can restrict our attention to black faces of size six.
The final charge of a black face $f$
of size six is at most one unless the face $f$ is isomorphic
to one of the faces depicted in Figure~\ref{fig-black}---note that
the amount of charge of $f$ can exceed $1$ unit only
if Rule A applies three times to $f$, Rule B applies twice to $f$ or
both Rules A and B apply to $f$. We analyze each of the configurations
separately in a series of three lemmas.

\begin{lemma}
\label{lm-black-AAA}
Let $C$ be a longest cycle of a fullerene graph $G$. Assume that
the discharging rules as described in Section~\ref{sect-rules}
have been applied. If $f=v_1v_2v_3v_4v_5v_6$ is a black face of $G$
such that the edges $v_5v_6$, $v_6v_1$, $v_1v_2$, $v_2v_3$ and $v_3v_4$
are contained in $C$ and the edge $v_4v_5$ is not,
then the final amount of charge of $f$ is at most $1$ unit.
\end{lemma}

\begin{figure}
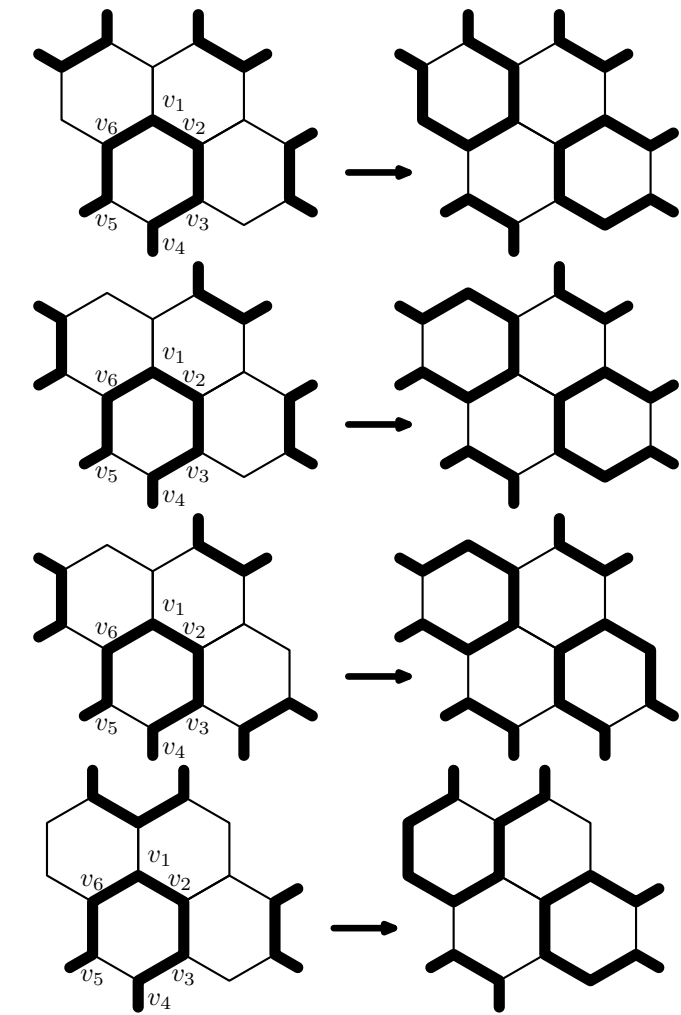

\begin{center}
\epsfbox{ful-ham.40}\\
\epsfbox{ful-ham.41}\\
\epsfbox{ful-ham.42}\\
\epsfbox{ful-ham.43}
\end{center}
\caption{Configurations analyzed in the proof of Lemma~\ref{lm-black-AAA}.}
\label{fig-black-AAA}
\end{figure}

\begin{proof}
The face $f$ can receive charge
only by Rule~A from the faces $f_{61}$, $f_{12}$ and $f_{23}$. Assume
for the sake of contradiction that $f$ receives charge of $1/2$ unit
from each of these three faces. In particular, $G$ contains, up to symmetry,
one of the configurations depicted in Figure~\ref{fig-black-AAA} (recall
that $G$ cannot contain a path formed by three white vertices
by Lemma~\ref{lm-path}).
Rerouting the cycle $C$ as indicated in the figure yields a cycle of $G$
longer than the cycle $C$ which contradicts our choice of $C$.
We conclude that Rule~A can apply at most twice to the face $f$.
\end{proof}

\begin{lemma}
\label{lm-black-AB}
Let $C$ be a longest cycle of a fullerene graph $G$. Assume that
the discharging rules as described in Section~\ref{sect-rules}
have been applied. If $f=v_1v_2v_3v_4v_5v_6$ is a black face of $G$
such that the edges $v_2v_3$, $v_4v_5$, $v_5v_6$ and $v_6v_1$
are contained in $C$ and the edges $v_1v_2$ and $v_3v_4$ are not,
then the final amount of charge of $f$ is at most $1$ unit.
\end{lemma}

\begin{figure}
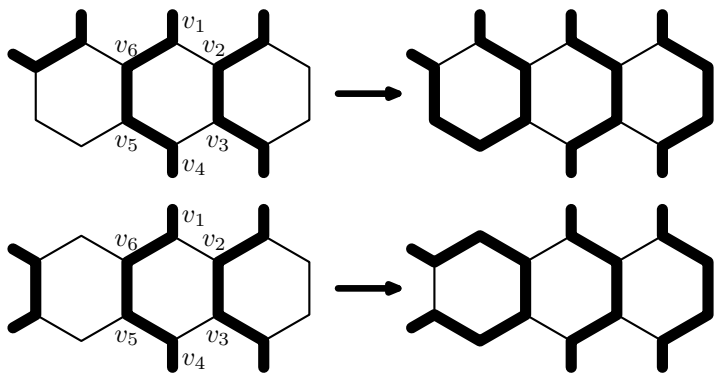

\begin{center}
\epsfbox{ful-ham.38}\\[3mm]
\epsfbox{ful-ham.39}
\end{center}
\caption{Configurations analyzed in the proof of Lemma~\ref{lm-black-AB}.}
\label{fig-black-AB}
\end{figure}

\begin{proof}
If the final amount of charge of $f$ is greater than $1$ unit,
then $f$ receives $1$ unit of charge from the face $f_{23}$ by Rule B and
$1/2$ unit of charge from the face $f_{56}$ by Rule A. Hence, $G$
contains one of the two configurations depicted in Figure~\ref{fig-black-AB}.
In either of the two cases, it is possible to reroute the cycle $C$ as
indicated in Figure~\ref{fig-black-AB} to obtain a cycle of $G$ longer than
$C$, a contradiction.
\end{proof}

\begin{lemma}
\label{lm-black-BB}
Let $C$ be a longest cycle of a fullerene graph $G$. Assume that
the discharging rules as described in Section~\ref{sect-rules}
have been applied. If $f=v_1v_2v_3v_4v_5v_6$ is a black face of $G$
such that the edges $v_2v_3$, $v_4v_5$ and $v_6v_1$ are contained in $C$ and
the edges $v_1v_2$, $v_3v_4$ and $v_5v_6$ are not,
then the final amount of charge of $f$ is at most $1$ unit.
\end{lemma}

\begin{figure}
\begin{center}
\epsfbox{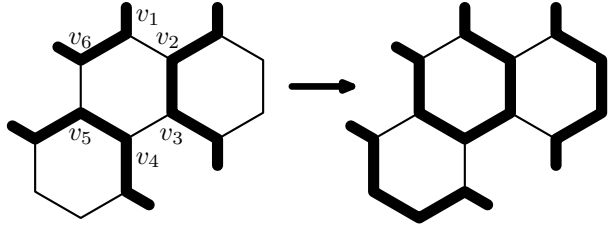}
\end{center}
\caption{The configuration analyzed in the proof of Lemma~\ref{lm-black-BB}.}
\label{fig-black-BB}
\end{figure}

\begin{proof}
The face $f$ can receive charge only by Rule~B. Assume that
Rule B applies twice to $f$. By symmetry, we may assume
that the charge is given by the faces $f_{23}$ and $f_{45}$.
In particular, the graph $G$ contains the configuration depicted
in Figure~\ref{fig-black-BB}. Reroute now the cycle $C$ as
indicated in the figure. Since the obtained cycle is longer than
the cycle $C$, we conclude that Rule~B cannot apply twice
to the face $f$.
\end{proof}

Lemmas~\ref{lm-black-AAA}, \ref{lm-black-AB} and~\ref{lm-black-BB}
yield the following.

\begin{lemma}
\label{lm-black}
Let $C$ be a longest cycle of a fullerene graph $G$. Assume that
the discharging rules as described in Section~\ref{sect-rules}
have been applied. The final amount of charge of any black face
of $G$ is at most $1$ unit.
\end{lemma}

%
%

\section{Main result}
\label{sect-main}

\begin{theorem}
\label{thm-main}
Let $G$ be a fullerene graph with $n$ vertices.
The graph $G$ contains a cycle of length at least $5n/6-2/3$.
\end{theorem}

\begin{proof}
Consider a longest cycle $C$ contained in the graph $G$ and
apply the discharging procedure described in Section~\ref{sect-rules}.
By Lemmas~\ref{lm-white} and~\ref{lm-black}, every white face and every
black face a has final charge of at most $1$ unit. Since the initial amount of
charge of other faces is $1$ unit and the other faces do not send out or
receive any charge, we conclude that the final amount of charge
of any face of $G$ is at most $1$ unit.

Each white vertex has initially been assigned $3$ units of charge.
Since the final amount of charge of every face is at most $1$ unit,
the amount of charge was preserved during the discharging phase and
vertices do not have any charge at the end of the process, there are
at most $f/3$ white vertices where $f$ is the number of faces of $G$.
By Euler's formula, $n=2f-4$. Hence, there are at most $n/6+2/3$ white
vertices. Consequently, there are at least $5n/6-2/3$ black vertices and
thus the length of the cycle $C$ is at least $5n/6-2/3$.
\end{proof}

\end{document}